\documentclass[12pt,twoside]{article}

\usepackage{latexsym}
\usepackage{amsmath}
\usepackage{amssymb}
\usepackage{pstricks}
\usepackage{pst-node}
\usepackage{multido}
\usepackage{calc}
\usepackage{enumerate}
\usepackage{subfigure}
\usepackage{color}
\usepackage{multirow}

\newtheorem{theorem}{Theorem}[section]

\newtheorem{conj}[theorem]{Conjecture}

\newcommand{\Prf}{\noindent{\bf Proof}.\quad }

\newcommand{\qed}{\hfill$\Box$}

\title{Irreducible pseudo $2$--factor isomorphic cubic bipartite graphs}
\author{{\rm M. Abreu,} \\
\small Dipartimento di Matematica, Universit\`a della
        Basilicata, \\
      \small C. da Macchia Romana, 85100 Potenza,
         Italy.\\
         \small e-mail: marien.abreu@unibas.it \\
\\
{\rm D. Labbate},\\
\small Dipartimento di Matematica, Politecnico di Bari\\
      \small I-70125 Bari, Italy.\\
        \small e-mail: labbate@poliba.it \\
\\
{\rm J. Sheehan,}\\
      \small Department of Mathematical Sciences, King's College,\\
      \small Old Aberdeen AB24 3UE,
         Scotland.\\
        \small e-mail: j.sheehan@maths.abdn.ac.uk }
\date{}

\newcommand{\grGT}{
\cnode*(0,4){2.5pt}{U1}
\cnode(1,4){2.5pt}{W1}
\cnode*(2,4){2.5pt}{X1}
\cnode(3,4){2.5pt}{Y1}
\cnode*(4,4){2.5pt}{Z1}
\cnode(5,4){2.5pt}{V1}

\cnode(2,3){2.5pt}{T1}
\cnode*(5,3){2.5pt}{T2}

\cnode*(0,2){2.5pt}{U2}
\cnode(1,2){2.5pt}{W2}
\cnode*(2.4,2){2.5pt}{X2}
\cnode(3,2){2.5pt}{Y2}
\cnode*(4,2){2.5pt}{Z2}
\cnode(5.5,2){2.5pt}{V2}

\cnode*(0,1){2.5pt}{U3}
\cnode(1,1){2.5pt}{W3}
\cnode*(2,1){2.5pt}{X3}
\cnode(3,1){2.5pt}{Y3}
\cnode*(4,1){2.5pt}{Z3}
\cnode(5,1){2.5pt}{V3}

\ncline{U1}{W1}
\ncline{W1}{X1}
\ncline{X1}{Y1}
\ncline{Y1}{Z1}
\ncline{Z1}{V1}

\ncline{U3}{W3}
\ncline{W3}{X3}
\ncline{X3}{Y3}
\ncline{Y3}{Z3}
\ncline{Z3}{V3}

\ncline{U1}{W2}
\ncline{U2}{W1}
\ncline{U2}{W3}
\ncline{U3}{W2}
\ncline{W2}{X2}
\ncline{X2}{Y2}

\ncline{Y1}{Z2}
\ncline{Y2}{Z1}
\ncline{Y2}{Z3}
\ncline{Y3}{Z2}
\ncline{Z2}{V2}

\ncline{X1}{T1}
\ncline{X2}{T1}
\ncline{X3}{T1}

\ncline{V1}{T2}
\ncline{V2}{T2}
\ncline{V3}{T2}
}

\newcommand{\grGTVert}{
\cnode*(0,4){2.5pt}{U1}
\cnode(1,4){2.5pt}{W1}
\cnode*(2,4){2.5pt}{X1}
\cnode(3,4){2.5pt}{Y1}
\cnode*(4,4){2.5pt}{Z1}
\cnode(5,4){2.5pt}{V1}

\cnode(2,3){2.5pt}{T1}
\cnode*(5,3){2.5pt}{T2}

\cnode*(0,2){2.5pt}{U2}
\cnode(1,2){2.5pt}{W2}
\cnode*(2.4,2){2.5pt}{X2}
\cnode(3,2){2.5pt}{Y2}
\cnode*(4,2){2.5pt}{Z2}
\cnode(5.5,2){2.5pt}{V2}

\cnode*(0,1){2.5pt}{U3}
\cnode(1,1){2.5pt}{W3}
\cnode*(2,1){2.5pt}{X3}
\cnode(3,1){2.5pt}{Y3}
\cnode*(4,1){2.5pt}{Z3}
\cnode(5,1){2.5pt}{V3}
}

\begin{document}
\maketitle

\begin{abstract}
A bipartite graph is {\em pseudo $2$--factor isomorphic} if all
its $2$--factors have the same parity of number of circuits.
In \cite{ADJLS} we proved that the only essentially
$4$--edge-connected pseudo $2$--factor isomorphic cubic bipartite
graph of girth $4$ is $K_{3,3}$, and conjectured \cite[Conjecture 3.6]{ADJLS} that the only
essentially $4$--edge-connected cubic bipartite graphs are
$K_{3,3}$, the Heawood graph and the Pappus graph.

There exists a characterization of symmetric configurations $n_3$
due to Martinetti (1886) in which all symmetric  configurations $n_3$ can be
obtained from an infinite set of so called {\em irreducible}
configurations \cite{VM}. The list of irreducible configurations
has been completed by Boben \cite{B} in terms of their {\em
irreducible Levi graphs}.

In this paper we characterize irreducible pseudo $2$--factor isomorphic cubic bipartite graphs
proving that the only pseudo $2$--factor isomorphic irreducible Levi graphs
are the Heawood and Pappus graphs. Moreover, the obtained characterization allows us to partially prove the above Conjecture.

\end{abstract}

\section{Introduction}

All graphs considered are finite and simple (without loops or
multiple edges). A graph with a $2$--factor is said to be
{\em $2$--factor hamiltonian} if all its $2$--factors are
Hamilton circuits, and, more generally,
{\em $2$--factor isomorphic} if all its
$2$--factors are isomorphic. Examples of such graphs are $K_4$,
$K_5$, $K_{3,3}$, the Heawood graph (which are all $2$--factor
hamiltonian) and the Petersen graph (which is $2$--factor
isomorphic). Several recent papers have addressed the problem of characterizing
families of graphs (particularly regular graphs) which have these
properties. It is shown in \cite{AFJLS,FJLS2} that $k$--regular
$2$--factor isomorphic bipartite graphs exist only when $k\in\{2,3\}$ and an
infinite family of $3$--regular $2$--factor hamiltonian bipartite
graphs, based on $K_{3,3}$ and the Heawood graph, is constructed
in \cite{FJLS2}. It is conjectured in \cite{FJLS2} that every
$3$--regular $2$--factor hamiltonian bipartite graph belongs to
this family. Faudree, Gould and Jacobsen in \cite{FGJ} determine
the maximum number of edges in both $2$--factor hamiltonian graphs
and $2$--factor hamiltonian bipartite graphs.
In addition, Diwan \cite{Di} has shown that $K_4$  is the only
$3$--regular $2$--factor hamiltonian planar graph.

Moreover, $2$--factor isomorphic bipartite graphs are extended in \cite{ADJLS}
to the more general family of {\em pseudo $2$--factor isomorphic graphs}
i.e. graphs $G$ with the property that the parity of the number of circuits in a
$2$--factor is the same for all $2$--factors of $G$. Example of these graphs are
$K_{3,3}$, the Heawood graph $H_0$ and the Pappus graph $P_0$.
Finally, it is proved in \cite{ADJLS} that pseudo $2$--factor isomorphic
$2k$--regular graphs and $k$--regular digraphs do not exist
for $k\geq 4$. Recently these results has been generalized in \cite{ALS}
for regular graphs which are not necessarily bipartite.

An incidence structure is {\em linear} if two different points are
incident with at most one line. A  {\em  symmetric configuration} $n_k$ (or {\em $n_k$ configuration})
is a linear incidence structure consisting of $n$ points
and $n$ lines such that each point and line is respectively
incident with $k$ lines and points.
Let $\cal C$ be a symmetric configuration $n_k$, its {\em Levi graph} $G({\cal C})$
is a $k$--regular bipartite graph whose vertex set are the points and the lines of $\cal C$ and
there is an edge between a point and a line in the graph if and only if they are incident in $\cal C$.
We will indistinctly refer to Levi graphs of configurations as their {\em incidence graphs}.

Let $G$ be a graph and $E_1$ be an edge-cut of $G$.
We say that $E_1$ is a {\em non-trivial edge-cut} if all
components of $G-E_1$ have at least two vertices.
The graph $G$ is {\em essentially $4$--edge--connected} if
$G$ is $3$--edge--connected and has no non-trivial $3$--edge--cuts.

\begin{conj}\label{conjADJLS} \cite{ADJLS}
Let $G$ be an essentially $4$--edge--connected pseudo $2$--factor
isomorphic cubic bipartite graph. Then $G \in \{K_{3,3},H_0,P_0\}.$
\end{conj}


\begin{theorem}\label{thmADJLS}\cite{ADJLS}
Let $G$ be an essentially $4$--edge--connected pseudo $2$--factor
isomorphic cubic bipartite graph. Suppose that $G$ contains a $4$--circuit.
Then $G = K_{3,3}.$ \qed
\end{theorem}

It follows from Theorem \ref{thmADJLS} that an essentially
$4$--edge--connected pseudo $2$--factor isomorphic cubic bipartite
graph of girth greater than or equal to $6$ is the Levi graph of
a symmetric configuration $n_3$. In $1886$ V. Martinetti \cite{VM}
characterized symmetric configurations $n_3$, showing that they can be
obtained from an infinite set of so called {\em irreducible}
configurations, of which he gave a list.
Recently, Boben proved that Martinetti's list of irreducible
configurations was incomplete and completed it \cite{B}.
Boben's list of irreducible configurations was obtained characterizing
their Levi graphs, which he called {\em irreducible Levi graphs} (cf. Section \ref{symconf}).

In this paper, we characterize {\em irreducible} pseudo $2$--factor isomorphic cubic bipartite graphs proving that
{\em the Heawood and the Pappus graphs are the only irreducible
Levi graphs which are pseudo $2$--factor isomorphic cubic bipartite.}
Moreover, the obtained characterization allows us to partially prove Conjecture \ref{conjADJLS}, i.e. in the case of irreducible pseudo $2$--factor isomorphic cubic bipartite graphs.

\section{Symmetric Configurations $\bf n_3$}\label{symconf}

In 1886, Martinetti \cite{VM} provided a construction for a symmetric configuration $n_3$ from a symmetric configuration $(n-1)_3$, say $\cal C$. Suppose
that in $\cal C$ there exist two parallel (non--intersecting)
lines $l_1 = \{\alpha, \alpha_1, \alpha_2\}$ and
$r_1 = \{\beta, \beta_1, \beta_2\}$ such that the points
$\alpha$ and $\beta$ are not on a common line. Then a symmetric configuration $n_3$ is obtained from $\cal C$ by deleting the lines
$l_1,r_1$, adding a point $\mu$ and adding the
lines $h_1 = \{\mu, \alpha_1, \alpha_2\}$,
$h_2 = \{\mu, \beta_1, \beta_2\}$ and $h_3 = \{\mu, \alpha, \beta\}$.
Not all symmetric configurations $n_3$ can be obtained using this method
on some symmetric configuration $(n-1)_3$. The configurations that cannot
be obtained in this way are called {\em irreducible
configurations}, while the others are {\em reducible configurations}.
However, if all irreducible symmetric configurations $n_3$ are known, then all
symmetric configurations $n_3$ can be constructed iteratively with Martinetti's
method. The list of irreducible configurations in \cite{VM} turned out to
be incomplete and it has been recently completed by Boben in \cite[Thm. 8]{B}.

\begin{theorem} \cite{B}
All connected irreducible $n_3$ configurations  are:

\begin{enumerate}
\item cyclic configurations with base line $\{0,1,3\}$;
\item the configurations with their incidence graphs
$T_1(n)$, $T_2(n)$, $T_3(n)$, \linebreak $n \geq 1$, each of them
giving precisely one $(10n)_3$ configuration, and;
\item the Pappus configuration.
\end{enumerate}
\end{theorem}

As mentioned before, Boben's list was obtained by studying the Levi graphs of
irreducible configurations, which are called {\em irreducible
Levi graphs}. Such graphs turned out to be either the Pappus graph,
or belong to one of four infinite families $D(n)$, $T_1(n)$,
$T_2(n)$, $T_3(n)$, $n \geq 1$, which we now proceed to
describe.

\

{\bf The $D(n)$ family:} Let $C(n)$, $n \geq 1$, be the graph on $6n$ vertices, consisting of $n$
segments ($6$--circuits labeled as in Fig. 1), linked by the edges
$v_{i-1}^1u_i^1$, $v_{i-1}^2u_i^4$, and $u_{i-1}^3u_i^2$, for $i \geq
2$.

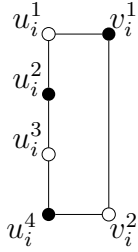
\begin{figure}[htb]
\centering
\begin{pspicture}(2,3.5) \psset{unit=0.8, linewidth=0.05pt}
\cnode(0,4){2.5pt}{N1}\rput(-0.35,4.25){$u_i^1$}
\cnode*(0,3){2.5pt}{N2}\rput(-0.35,3.25){$u_i^2$}
\cnode(0,2){2.5pt}{N3}\rput(-0.35,2.25){$u_i^3$}
\cnode*(0,1){2.5pt}{N4}\rput(-0.45,0.75){$u_i^4$}

\cnode*(1,4){2.5pt}{M1}\rput(1.25,4.25){$v_i^1$}
\cnode(1,1){2.5pt}{M2}\rput(1.25,0.75){$v_i^2$}

\ncline{N1}{N2}
\ncline{N2}{N3}
\ncline{N3}{N4}

\ncline{N1}{M1}
\ncline{N4}{M2}

\ncline{M1}{M2}
\end{pspicture}
\caption{label for the $6$--circuits of the definition of $C(n)$}
\end{figure}

Let the graph $D(n)$, $n \geq 1$ be defined as follows:

For $n \equiv 0 \, (mod \, 3)$ let $D(n)$ be the graph $C(m)$,
$m=n/3$, with the edges $u_1^1v_m^1$, $u_1^4v_m^2$, and
$u_1^2u_m^3$ added (cf. Fig. 2).

For $n \equiv 1 \, (mod \, 3)$ let $D(n)$ be the graph $C(m)$,
$m=(n-1)/3$, with two vertices $w_m^1,w_m^2$ and the edges $u_1^1w_m^1$,
$u_1^2v_m^2$, $u_1^4w_m^2$, $w_m^1w_m^2$, $w_m^1u_m^3$,
$w_m^2v_m^1$ added.

For $n \equiv 2 \, (mod \, 3)$ let $D(n)$ be the graph $C(m)$,
$m=(n-2)/3$, with four vertices $w_m^1,w_m^2,w_m^3,w_m^4$ and the edges
$v_m^1w_m^1$, $v_m^2w_m^4$, $u_m^3w_m^2$, $u_1^1w_m^4$,
$u_1^2w_m^1$, $u_1^4w_m^3$, $w_m^1w_m^2$, $w_m^2w_m^3$, $w_m^3w_m^4$ added.

\begin{figure}[htb]
\centering
\begin{pspicture}(4,4) \psset{unit=0.8, linewidth=0.05pt}
\cnode(0,4){2.5pt}{N1}\rput(-0.45,4.25){$u_1^1$}
\cnode*(1,4){2.5pt}{N2}\rput(0.65,4.25){$v_1^1$}
\cnode(2,4){2.5pt}{N3}\rput(1.65,4.25){$u_2^1$}
\cnode*(3,4){2.5pt}{N4}\rput(2.65,4.25){$v_2^1$}
\cnode(4,4){2.5pt}{N5}\rput(3.65,4.25){$u_3^1$}
\cnode*(5,4){2.5pt}{N6}\rput(4.65,4.25){$v_3^1$}

\cnode*(0,3){2.5pt}{N7}\rput(-0.45,3.25){$u_1^2$}
\cnode*(2,3){2.5pt}{N8}\rput(1.65,3.25){$u_2^2$}
\cnode*(4,3){2.5pt}{N9}\rput(3.65,3.25){$u_3^2$}

\cnode(0,2){2.5pt}{N10}\rput(-0.45,2.25){$u_1^3$}
\cnode(2,2){2.5pt}{N11}\rput(1.65,2.25){$u_2^3$}
\cnode(4,2){2.5pt}{N12}\rput(3.70,2.25){$u_3^3$}

\cnode*(0,1){2.5pt}{N13}\rput(-0.45,0.75){$u_1^4$}
\cnode(1,1){2.5pt}{N14}\rput(0.65,0.75){$v_1^2$}
\cnode*(2,1){2.5pt}{N15}\rput(1.65,0.75){$u_2^4$}
\cnode(3,1){2.5pt}{N16}\rput(2.65,0.75){$v_2^2$}
\cnode*(4,1){2.5pt}{N17}\rput(3.65,0.75){$u_3^4$}
\cnode(5,1){2.5pt}{N18}\rput(4.65,0.75){$v_3^2$}

\ncline{N1}{N2}
\ncline{N1}{N7}

\nccurve[angleA=120,angleB=60]{N1}{N6}

\ncline{N2}{N3}
\ncline{N2}{N14}

\ncline{N3}{N4}
\ncline{N3}{N8}

\ncline{N4}{N5}
\ncline{N4}{N16}

\ncline{N5}{N6}
\ncline{N5}{N9}

\ncline{N6}{N18}

\ncline{N7}{N10}
\ncline{N7}{N12}

\ncline{N8}{N10}
\ncline{N8}{N11}

\ncline{N9}{N11}
\ncline{N9}{N12}

\ncline{N10}{N13}

\ncline{N11}{N15}

\ncline{N12}{N17}

\ncline{N13}{N14}

\nccurve[angleA=-120,angleB=-60]{N13}{N18}

\ncline{N14}{N15}

\ncline{N15}{N16}

\ncline{N16}{N17}

\ncline{N17}{N18}

\end{pspicture}
\caption{The graph $D(9)$}
\end{figure}
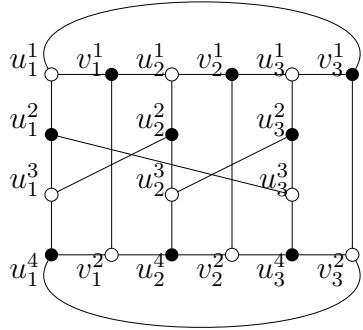

A {\em cyclic configuration} has ${\mathbb Z}_n=\{0,1, \ldots, n-1\}$ as set of
points and ${\cal B}= \{ \{0,b,c\}, \{1,b+1,c+1\}, \ldots, \{n-1,b+n-1,c+n-1\}
\}$ as set of lines, where the operations are modulo $n$, and the
{\em base line} is $\{0,b,c\}$ for $b,c \in {\mathbb Z}_n$.

Note that the graphs $D(n)$ are the Levi graphs of the cyclic $n_3$
configurations with base line $\{0,1,3\}$. In particular,
for $n=7$ the cyclic $7_3$ configuration is the Fano plane and $D(7)$ is the Heawood
graph $H_0$.

\

{\bf The $T_1(n)$, $T_2(n)$ and $T_3(n)$ families:} Let $T(n)$, $n \geq 1$, be the graph on $20n$ vertices consisting
of $n$ segments $G_T$ shown in Fig. 3, linked by the edges
$v_{i-1}^1u_i^1$, $v_{i-1}^2u_i^2$, $v_{i-1}^3u_i^3$, for $i \geq
2$.

\begin{figure}[htb]
\centering
\begin{pspicture}(3.5,3.5) \psset{unit=0.8, linewidth=0.05pt}
\grGT{}
\rput(-0.35,4){$u_i^1$}
\rput(1.25,4.35){$w_i^1$}
\rput(2.25,4.35){$x_i^1$}
\rput(3.25,4.35){$y_i^1$}
\rput(4.25,4.35){$z_i^1$}
\rput(5.5,4){$v_i^1$}
\rput(2.25,3.25){$t_i^1$}
\rput(5.25,3.25){$t_i^2$}
\rput(-0.35,2){$u_i^2$}
\rput(1.25,2.25){$w_i^2$}
\rput(2.5,1.65){$x_i^2$}
\rput(2.95,2.45){$y_i^2$}
\rput(4.25,2.25){$z_i^2$}
\rput(5.9,2){$v_i^2$}
\rput(-0.35,1){$u_i^3$}
\rput(1.25,0.75){$w_i^3$}
\rput(2.25,0.75){$x_i^3$}
\rput(3.25,0.75){$y_i^3$}
\rput(4.25,0.75){$z_i^3$}
\rput(5.5,1){$v_i^3$}
\end{pspicture}
\caption{Segment graph $G_T$}
\end{figure}
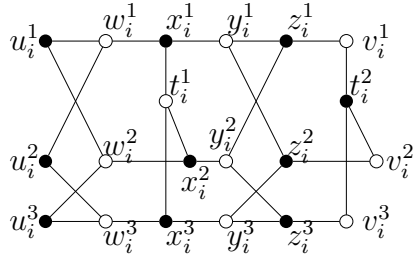

Let $T_1(n)$ be the graph obtained from $T(n)$ by adding the
edges $u_1^1v_n^1$, $u_1^2v_n^2$, $u_1^3v_n^3$.
Let $T_2(n)$ be the graph obtained from $T(n)$ by adding the
edges $u_1^3v_n^1$, $u_1^2v_n^2$, $u_1^1v_n^3$.
Let $T_3(n)$ be the graph obtained from $T(n)$ by adding the
edges $u_1^1v_n^3$, $u_1^2v_n^1$, $u_1^3v_n^2$.
In \cite{B}, Boben proved that for each fixed value of $n$, no two
of the graphs $T_1(n)$, $T_2(n)$, $T_3(n)$ are isomorphic.

\

Note that $T_1(1)$ is the Levi graph of Desargues' configuration,
and $T_2(1)$, $T_3(1)$ correspond to the Levi graphs of the
configurations $10_3F$ and $10_3G$ 
respectively according to Kantor's \cite{Kantor} notation for the ten
$10_3$ configurations.

\

{\bf The Pappus graph:} Recall that the Levi graph of the Pappus $9_3$ configuration
is the following pseudo $2$--factor isomorphic but not $2$--factor isomorphic cubic bipartite graph \cite{ADJLS}, called the {\em Pappus graph $P_0$}.

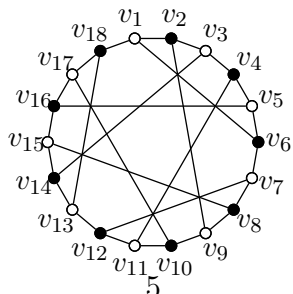
\begin{figure}[h]
\centering
\begin{pspicture}(0,1.3) \SpecialCoor \psset{unit=1.4}
            \degrees[-18]
            \multido{\ia=13+2,\ib=12+2}{18}{%
            \cnode[linewidth=0.5pt](1;\ia){2.3pt}{M\ia}%
            \cnode*(1;\ib){2.3pt}{M\ib}}

            \multido{\ic=13+1,\id=1+1,\ie=12+1}{18}{%
            \rput(1.20;\ic){{\small$v_{\id}$}}
            \ncline[linewidth=0.5pt]{M\ic}{M\ie} }

            \ncline[linewidth=0.5pt]{M13}{M18}
            \ncline[linewidth=0.5pt]{M14}{M21}
            \ncline[linewidth=0.5pt]{M15}{M26}
            \ncline[linewidth=0.5pt]{M16}{M23}
            \ncline[linewidth=0.5pt]{M17}{M28}
            \ncline[linewidth=0.5pt]{M19}{M24}
            \ncline[linewidth=0.5pt]{M20}{M27}
            \ncline[linewidth=0.5pt]{M22}{M29}
            \ncline[linewidth=0.5pt]{M25}{M30}
\end{pspicture}
\end{figure}


\section{$2$--factors of Irreducible Levi Graphs}


Let $G$ be a graph and $u,v$ be two vertices in $G$.
Then a {\em $(u,v)$--path} is a path from $u$ to $v$.
Given two disjoint paths $P=u_1, \ldots, u_n$
and $Q=u_{n+1}, \ldots, u_{n+m}$ (except maybe for $u_1=u_{n+m}$),
the path $PQ=u_1, \ldots, u_{n+m}$ is the concatenation of $P$
and $Q$ together with the edge $u_nu_{n+1}$.
Similarly, for a vertex $v \in (G - P) \cup \{v_1\}$, the path
$Pv$ is composed by $P$, $v$ and the edge $u_nv$
If $u_1=u_{n+m}$ or $u_1=v$ we write $(PQ)$ and $(Pv)$ respectively,
to emphasize that $PQ$ and $Pv$ are circuits.

    \begin{theorem}\label{MainThm}
The Heawood and the Pappus graphs are the only irreducible
Levi graphs which are pseudo $2$--factor isomorphic.
    \end{theorem}

\Prf  It is straightforward to show that the Heawood graph $H_0$ is $2$--
factor hamiltonian and hence pseudo 2-factor isomorphic (cf. \cite{FJLS2}).
We have already proved in \cite[Proposition 3.3]{ADJLS} that the Pappus
graph is pseudo $2$--factor isomorphic. We need to prove
that all other irreducible Levi graphs are not pseudo $2$--factor isomorphic
and we will do so by finding two $2$--factors with different parity of number of circuits
in each of them.

\noindent The following paths will be used for constructing $2$--factors in
$D(n)$, for $n \geq 8$.

\small{
$$
\begin{array}{ll}
L_i^1=u_i^1u_i^2u_i^3u_i^4v_i^2v_i^1                             & L_i^2=u_i^4u_i^3u_i^2u_i^1v_i^1v_i^2 \\
M_1=u_1^4v_1^2u_2^4u_2^3u_3^2u_3^3u_3^4v_2^2v_2^1u_3^1v_3^1v_3^2 & M_2=u_m^4v_m^2v_m^1u_m^1u_m^2u_m^3w_m^1w_m^2u_1^4 \\
N_i=u_i^2u_i^1v_{i-1}^1v_{i-1}^2u_i^4u_i^3                       & N_m=w_m^2w_m^1v_m^1v_m^2w_m^4w_m^3u_1^4v_1^2u_2^4u_2^3 \\
C_1=(u_1^1u_1^2u_1^3u_2^2u_2^1v_1^1u_1^1)                        & C_2=(u_1^1u_1^2u_1^3u_1^4w_m^2w_m^1u_1^1) \\
C_3=(v_1^1v_1^2u_2^4v_2^2v_2^1u_3^1v_3^1v_3^2u_3^4u_3^3u_3^2u_2^3u_2^2u_2^1v_1^1)
\end{array}
$$
}

\noindent Hamiltonian $2$--factors in $D(n)$ are
\small{ $ \left\{ \begin{array}{lc}
(L_1^1L_2^1 \cdots L_m^1u_1^1)                     & n \equiv 0 \mod \, 3 \\
(L_1^1L_2^1 \cdots L_m^1w_m^2w_m^1u_1^1)           & n \equiv 1 \mod \, 3 \\
(L_1^1L_2^1 \cdots L_m^1w_m^1w_m^2w_m^3w_m^4u_1^1) & n \equiv 2 \mod \, 3 \\
\end{array}
\right.$ }

%

\noindent Disconnected $2$--factors with exactly two circuits in $D(n)$ are

\small{ $$ \left\{ \begin{array}{llc}
C_1 \cup (M_1u_1^4)                        & n = 9                  & \multirow{2}*{\quad for \, $n \equiv 0 \mod 3$}\\
C_1 \cup (M_1L_4^2 \cdots L_m^2u_1^4)      & n =3m, \, m \ge 4      \\
\\
C_2 \cup C_3                               & n = 10                 & \multirow{3}*{\quad for \, $n \equiv 1 \mod 3$}\\
C_1 \cup (M_1M_2)                          & n = 13                 \\
C_1 \cup (M_1L_4^2 \cdots L_{m-1}^2M_2)    & n = 3m + 1, \, m \ge 5 \\
\\
C_1 \cup (N_mw_m^2)                        & n = 8                  & \multirow{2}*{\quad for \, $n \equiv 2 \mod 3$}\\
C_1 \cup (N_3 \cdots N_mu_3^2)             & n = 3m + 2, \, m \ge 3 \\
\end{array}
\right.$$}

%
%

\

\noindent Now we need to find such pairs of $2$--factors for the graphs
$T_1(n)$, $T_2(n)$ and $T_1(n)$, $n \geq 1$. To this purpose
we need to consider the following six paths in the segment graph
$G_T$ from Fig 3.

\

{\small $P_i^1=u_i^2,w_i^1,u_i^1,w_i^2,u_i^3,w_i^3,x_i^3,y_i^3,z_i^3,
v_i^3,t_i^2,v_i^1,z_i^1,y_i^2,x_i^2,t_i^1,x_i^1,y_i^1,z_i^2,v_i^2$}

{\small $P_i^2=u_i^2,w_i^1,u_i^1,w_i^2,u_i^3,w_i^3,x_i^3,y_i^3,z_i^2,v_i^2$}

{\small $(P_i^3)=(v_i^1,t_i^2,v_i^3,z_i^3,y_i^2,x_i^2,t_i^1,x_i^1,y_i^1,z_i^1,v_i^1)$}

{\small $Q_i^1=u_i^3,w_i^3,u_i^2,w_i^1,u_i^1,w_i^2,x_i^2,y_i^2,z_i^3,
v_i^3,t_i^2,v_i^1,z_i^1,y_i^1,x_i^1,t_i^1,x_i^3,y_i^3,z_i^2,v_i^2$}

{\small $Q_i^2=u_i^3,w_i^2,u_i^1,w_i^1,u_i^2,w_i^3,x_i^3,y_i^3,z_i^2,v_i^2$}

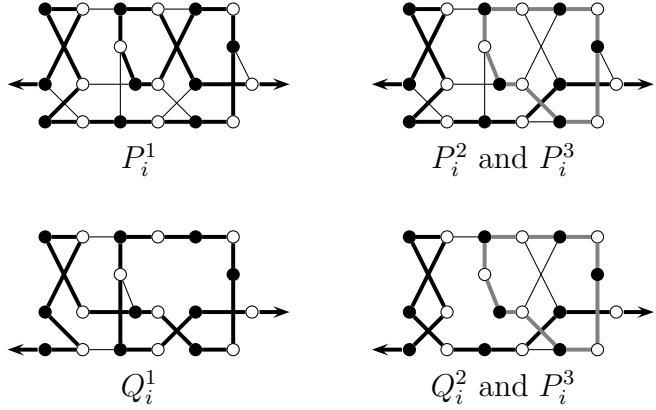
\rule{0cm}{2cm}
\begin{figure}[h]
\centering
\begin{tabular}{ccc}
    \subfigure{\label{fig5-a}\begin{pspicture}(3,-2.5) \psset{unit=0.5}
            \psset{linewidth=0.4pt}
            \grGTVert{}
            \ncline{W1}{X1} \ncline{Y1}{Z1} \ncline{U2}{W3}
            \ncline{W2}{X2} \ncline{Y2}{Z3} \ncline{Y3}{Z2}
            \ncline{X3}{T1} \ncline{V2}{T2}
            \psset{linewidth=1.5pt,linestyle=solid}
            \ncline{U2}{W1}\ncline{W1}{U1}\ncline{U1}{W2}
            \ncline{W2}{U3}\ncline{U3}{W3}\ncline{W3}{X3}
            \ncline{X3}{Y3}\ncline{Y3}{Z3}\ncline{Z3}{V3}
            \ncline{V3}{T2}\ncline{T2}{V1}\ncline{V1}{Z1}
            \ncline{Z1}{Y2}\ncline{Y2}{X2}\ncline{X2}{T1}
            \ncline{T1}{X1}\ncline{X1}{Y1}\ncline{Y1}{Z2}
            \ncline{Z2}{V2}
            \psline{->}(5.7,2)(6.5,2)\psline{->}(-0.2,2)(-1,2)
            \rput(2.5,0){$P_i^1$}
    \end{pspicture}} & \hspace{1cm} &%
    \subfigure{\label{fig5-b}\begin{pspicture}(3,-2.5) \psset{unit=0.5}
            \psset{linewidth=0.4pt}
            \grGTVert{}
            \ncline{Y3}{Z3} \ncline{Z1}{Y2} \ncline{Y1}{Z2}
            \ncline{W1}{X1} \ncline{U2}{W3} \ncline{W2}{X2}
            \ncline{X3}{T1} \ncline{V2}{T2}
            \psset{linewidth=1.5pt,linestyle=solid}
            \ncline{U2}{W1} \ncline{W1}{U1} \ncline{U1}{W2}
            \ncline{W2}{U3} \ncline{U3}{W3} \ncline{W3}{X3}
            \ncline{X3}{Y3} \ncline{Y3}{Z2} \ncline{Z2}{V2}
            \psline{->}(5.7,2)(6.5,2)\psline{->}(-0.2,2)(-1,2)
            \psset{linewidth=1.5pt,linestyle=solid, linecolor=gray}
            \ncline{V1}{T2} \ncline{T2}{V3} \ncline{V3}{Z3}
            \ncline{Z3}{Y2} \ncline{Y2}{X2} \ncline{X2}{T1}
            \ncline{T1}{X1} \ncline{X1}{Y1} \ncline{Y1}{Z1}
            \ncline{Z1}{V1}
            \rput(2.5,0){$P_i^2$ and $P_i^3$}
    \end{pspicture}} \\%
    \subfigure{\label{fig5-d}\begin{pspicture}(3,0) \psset{unit=0.5}
            \psset{linewidth=0.4pt}
            \grGTVert{}
            \ncline{U1}{W2} \ncline{W3}{X3} \ncline{V3}{T2}
            \ncline{X2}{T1} \ncline{W1}{X1} \ncline{Y1}{Z1}
            \ncline{Y2}{Z3} \ncline{Y3}{Z2}
            \psset{linewidth=1.5pt,linestyle=solid}
            \ncline{U3}{W3} \ncline{W3}{U2} \ncline{U2}{W1}
            \ncline{W1}{U1} \ncline{U1}{W2} \ncline{W2}{X2}
            \ncline{X2}{Y2} \ncline{Y2}{Z3} \ncline{Z3}{V3}
            \ncline{V3}{T2} \ncline{T2}{V1} \ncline{V1}{Z1}
            \ncline{Z1}{Y1} \ncline{Y1}{X1} \ncline{X1}{T1}
            \ncline{T1}{X3} \ncline{X3}{Y3} \ncline{Y3}{Z2}
            \ncline{Z2}{V2}
            \psline{->}(5.7,2)(6.5,2)\psline{->}(-0.2,1)(-1,1)
            \rput(2.5,0){$Q_i^1$}
    \end{pspicture}} & \hspace{1cm} &%
    \subfigure{\label{fig5-e}\begin{pspicture}(3,0) \psset{unit=0.5}
            \psset{linewidth=0.4pt}
            \grGTVert{}
            \ncline{Y2}{Z1} \ncline{Z2}{Y1} \ncline{Y3}{Z3}
            \ncline{U1}{W2} \ncline{W3}{X3} \ncline{V3}{T2}
            \ncline{X2}{T1} \ncline{W1}{X1}
            \psset{linewidth=1.5pt,linestyle=solid}
            \ncline{U3}{W2} \ncline{W2}{U1} \ncline{U1}{W1}
            \ncline{W1}{U2} \ncline{U2}{W3} \ncline{W3}{X3}
            \ncline{X3}{Y3} \ncline{Y3}{Z2} \ncline{Z2}{V2}
            \psline{->}(5.7,2)(6.5,2)\psline{->}(-0.2,1)(-1,1)
            \psset{linewidth=1.5pt,linestyle=solid, linecolor=gray}
            \ncline{V1}{T2} \ncline{T2}{V3} \ncline{V3}{Z3}
            \ncline{Z3}{Y2} \ncline{Y2}{X2} \ncline{X2}{T1}
            \ncline{T1}{X1} \ncline{X1}{Y1} \ncline{Y1}{Z1}
            \ncline{Z1}{V1}
            \rput(2.5,0){$Q_i^2$ and $P_i^3$}
    \end{pspicture}}
\end{tabular}
\end{figure}

\noindent The paths $P_i^1$ and $Q_i^1$ are hamiltonian
$(u_i^2,v_i^2)$ and $(u_i^3,v_i^2)$--paths,
respectively. The paths $P_i^2$ and $Q_i^2$ are
$(u_i^2,v_i^2)$, and $(u_i^3,v_i^2)$--paths on $10$ vertices,
respectively. Finally, $(P_i^3)$ is a $10$--circuit
in $G_T  - P_i^2=G_T  - Q_i^2$.

In $T_1(n)$ and $T_2(n)$ the hamiltonian $2$--factor
$F_1(n)=(P_1^1P_2^1 \cdots P_n^1u_1^2)$ and the disconeccted $2$--factor with exactly two circuits
$F_2(n)=(P_1^2P_2^1 \cdots P_n^1u_1^2) \cup (P_1^3)$ (even)
show that these graphs are not pseudo $2$--factor isomorphic.

Similarly, in $T_3(n)$ the hamiltonian $2$--factor
$F'_1(1)=(Q_1^1P_2^1 \cdots P_n^1u_1^3)$ and the disconnected $2$--factor with exactly two circuits
$F'_2(1)=(Q_1^2P_2^1 \cdots P_n^1u_1^3) \cup (P_1^3)$
show that these graphs are not pseudo $2$--factor isomorphic.
\qed



\

Note that Theorem \ref{MainThm} proves Conjecture \ref{conjADJLS} in the case of irreducible pseudo $2$--factor isomorphic
cubic bipartite graphs. In the next Section we show that Conjecture \ref{conjADJLS} cannot be directly proved
from Theorem \ref{MainThm} by extending it to reducible Levi graphs.

\section{$\bf 2$--factors in extensions and reductions of Levi graphs of $\bf n_3$ configurations}

Recall that a {\em Martinetti extension} can be described in
terms of graphs as follows:

Let $G_1$ be the Levi graph of a symmetric configuration $n_3$ and suppose that in $G_1$ there are two edges $e_1=x_1y_1$
and $e_2=x_2y_2$ with no common neighbours, then the
graph $G:=G_1 - \{e_1,e_2\} + \{u,v\} + \{ux_1,ux_2,vy_1,vy_2\}$
is the Levi graph of an  $(n+1)_3$ configuration.

\begin{figure}[htb]
\begin{pspicture}(-5,2)
\psset{linewidth=0.5pt, unit=0.5}
\cnode*(0,0){2.5pt}{N1}\rput(-0.75,0){$y_1$}%
\cnode*(2,0){2.5pt}{N2}\rput(2.75,0){$y_2$}%
\cnode(0,4){2.5pt}{N3}\rput(-0.75,4){$x_1$}%
\cnode(2,4){2.5pt}{N4}\rput(2.75,4){$x_2$}%
\ncline{N1}{N3}
\ncline{N2}{N4}

\rput(4,2){$\longrightarrow$}

\cnode*(6,0){2.5pt}{M1}\rput(5.25,0){$y_1$}%
\cnode*(8,0){2.5pt}{M2}\rput(8.75,0){$y_2$}%
\cnode(7,1){2.5pt}{M3}\rput(6.5,1){$v$}%
\cnode*(7,3){2.5pt}{M4}\rput(6.5,3){$u$}%
\cnode(6,4){2.5pt}{M5}\rput(5.25,4){$x_1$}%
\cnode(8,4){2.5pt}{M6}\rput(8.75,4){$x_2$}%
\ncline{M1}{M3}
\ncline{M2}{M3}
\ncline{M4}{M3}
\ncline{M4}{M5}
\ncline{M4}{M6}

\end{pspicture}
\caption{Martinetti Extension}\label{fig4}
\end{figure}
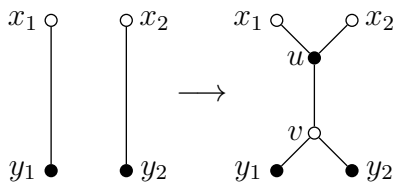

Similarly the Levi graph $G$ of a
symmetric configuration $(n+1)_3$ is {\em Martinetti reducible}
if there is an edge $e=uv$ in $G$ such that
either $G:=G_1 - \{u,v\} + x_1y_1 + x_2y_2$ or
$G:=G_1 - \{u,v\} + x_1y_2 + x_2y_1$
is again the Levi graph of a symmetric configuration $n_3$,
where $x_1,x_2,y_1,y_2$ are the neighbours of
$u$ and $v$ as in the following figure:

\begin{figure}[htb]
\begin{pspicture}(-3,2)
\psset{linewidth=0.5pt, unit=0.5}
\cnode*(0,0){2.5pt}{M1}\rput(-0.75,0){$y_1$}%
\cnode*(2,0){2.5pt}{M2}\rput(2.75,0){$y_2$}%
\cnode(1,1){2.5pt}{M3}\rput(0.5,1){$v$}%
\cnode*(1,3){2.5pt}{M4}\rput(0.5,3){$u$}%
\cnode(0,4){2.5pt}{M5}\rput(-0.75,4){$x_1$}%
\cnode(2,4){2.5pt}{M6}\rput(2.75,4){$x_2$}%

\ncline{M1}{M3}
\ncline{M2}{M3}
\ncline{M4}{M3}\Aput{$e$}
\ncline{M4}{M5}
\ncline{M4}{M6}

\rput(6,2){$\longrightarrow$}

\cnode*(8,0){2.5pt}{N1}\rput(7.25,0){$y_1$}%
\cnode*(10,0){2.5pt}{N2}\rput(10.75,0){$y_2$}%
\cnode(8,4){2.5pt}{N3}\rput(7.25,4){$x_1$}%
\cnode(10,4){2.5pt}{N4}\rput(10.75,4){$x_2$}%
\ncline{N1}{N3}
\ncline{N2}{N4}

\rput(12,2){$or$}

\cnode*(14,0){2.5pt}{P1}\rput(13.25,0){$y_1$}%
\cnode*(16,0){2.5pt}{P2}\rput(16.75,0){$y_2$}%
\cnode(14,4){2.5pt}{P3}\rput(13.25,4){$x_1$}%
\cnode(16,4){2.5pt}{P4}\rput(16.75,4){$x_2$}%
\ncline{P1}{P4}
\ncline{P2}{P3}

\end{pspicture}
\caption{Martinetti Reduction}\label{fig5}
\end{figure}
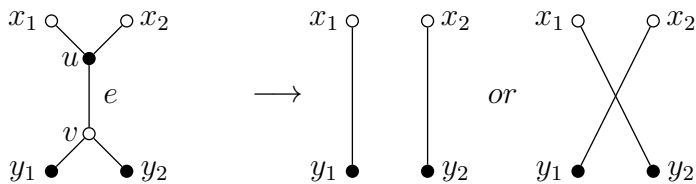

\

It is well known that the $7_3$ configuration, whose Levi graph is the Heawood graph, is not Martinetti extendible
and that the Pappus configuration is Martinetti extendible in a unique way; it is easy to show that this extension is not
pseudo $2$--factor isomorphic.

Let $\cal C$ be a symmetric configuration $n_3$ and
$\cal C'$ be a symmetric configuration $(n+1)_3$ obtained
from $\cal C$ through a Martinetti extension.
It can be easily checked that there are $2$--factors in $\cal C'$
that cannot be reduced to a $2$--factor in $\cal C$.
For example if $\cal C$ corresponds to the first option in Fig. \ref{fig5},
a $2$--factor of $\cal C'$ containing the path $x_1uvy_2$ will not reduce
to a $2$--factor in $\cal C$.
Conversely, there might be $2$--factors of $\cal C$ for
which the parity of number of circuits is not preserved
when extended to a $2$--factor in $\cal C'$.
For example, the graph $H_0*H_0$ (the star product \cite[p. 90]{HS} of the Heawood graph with itself)
which is $2$--factor hamiltonian and Martinetti
reducible (only through the edges of the non--trival
$3$--edge--cut), has all Martinetti reductions which are no longer pseudo
$2$--factor isomorphic.
Hence, we cannot directly prove Conjecture \ref{conjADJLS} by studying
the $2$--factors of reducible configurations from the
set of $2$--factors of their underlying irreducible ones.

\end{document}